\documentclass{amsart}

\newtheorem{theorem}{Theorem}
\newcommand{\bt}{\begin{theorem}}
\newcommand{\et}{\end{theorem}}

\newtheorem{lemma}{Lemma}
\newcommand{\bl}{\begin{lemma}}
\newcommand{\el}{\end{lemma}}
\newcommand{\beq}{\begin{equation}}
\newcommand{\eeq}{\end{equation}}
\newcommand{\benum}{\begin{enumerate}}
\newcommand{\eenum}{\end{enumerate}}

\newcommand{\mbc}{\ensuremath{\mathbf c}}
\newcommand{\mbr}{\ensuremath{\mathbf r}}

\newcommand{\mbz}{\ensuremath{\mathbf z}}
\newcommand{\mbo}{\ensuremath{\mathbf 0}}

\newcommand{\mbj}{\ensuremath{\mathbf j}}
\newcommand{\mcc}{\ensuremath{ \mathcal C}}
\newcommand{\mcr}{\ensuremath{ \mathcal R}}

\newcommand{\Rm}{\ensuremath{ \mathbf{R}^{m} }}
\newcommand{\Rn}{\ensuremath{ \mathbf{R}^n }}

\DeclareMathOperator{\kernel}{\text{kernel}}

\DeclareMathOperator{\diag}{\text{diag}}
\DeclareMathOperator{\image}{\text{image}}
\DeclareMathOperator{\qqand}{\qquad\text{and}\qquad}
\DeclareMathOperator{\rowsum}{\text{rowsum}}
\DeclareMathOperator{\colsum}{\text{colsum}}

\newcommand{\bmat}{\left(\begin{matrix}}
\newcommand{\emat}{\end{matrix}\right)}
\newcommand{\bsmallmat}{\left(\begin{smallmatrix}}
\newcommand{\esmallmat}{\end{smallmatrix}\right)}

\usepackage{amsmath,amssymb}\usepackage[all]{xy}\usepackage{tikz}

\date{\today}

\title{Sinkhorn limits in finitely many steps}
\author{Alex Cohen and Melvyn B. Nathanson}
\address{Yale University}
\email{alex.cohen@yale.edu}
\address{Lehman College (CUNY)}
\email{melvyn.nathanson@lehman.cuny.edu}

\subjclass[2010]{11C20, 11B75, 15B51, 05B20.}

\keywords{Matrix scaling, Sinkhorn limits, doubly stochastic matrices.}

\thanks{M.B.N. supported in part by a grant from the PSC-CUNY Research Award Program.}

\date{\today}
\begin{document}

\begin{abstract}
Applied to a nonnegative $m\times n$ matrix with a nonzero $\sigma$-diagonal, 
the sequence of matrices constructed by alternate row and column scaling conveges 
to a doubly stochastic matrix.  It is proved that if this sequence converges after only a finite number of scalings, 
then it converges after at most two scalings.
\end{abstract}

\maketitle

\section{The Sinkhorn-Knopp algorithm}

Let $A = \bmat a_{i,j} \emat$ be an $m\times n$ matrix.
We define the \emph{row sum}
\[
\rowsum_i(A) =  \sum_{j=1}^n a_{i,j} \qquad\text{for $i \in \{1,\ldots, m\}$}
\]
and the \emph{column sum}  
\[
\colsum_j(A) = \sum_{i=1}^m a_{i,j}  \qquad\text{for $j \in \{1,\ldots, n$\}.}
\]
Let $\mathbf{j_n} = \bsmallmat 1 \\ \vdots \\ 1 \esmallmat \in \Rn$ 
be the vector with all coordinates equal to 1.  
Let $A^t$ be the transpose of the matrix $A$.  
We have 
\[
A \mbj_n = \bmat \rowsum_1(A) \\ \vdots \\ \rowsum_m(A) \emat \in \Rm
\]
and
\[
A^t \mbj_m =  \bmat \colsum_1(A) \\ \vdots \\ \colsum_n(A) \emat \in \Rn.
\]

The $m\times n$ matrix $A = \bmat a_{i,j} \emat$ is \emph{nonnegative} if $a_{i,j} \geq 0$ 
for all $i \in \{1,\ldots, m\}$ and $j \in \{1,\ldots, n\}$.
The nonnegative $n \times n$ matrix $A$ is \emph{row stochastic} if 
$\rowsum_i(A) = 1$ for all $i \in \{1,\ldots, n\}$, 
and \emph{column stochastic} if 
$\colsum_j(A) = 1$  for all $ j\in \{1,\ldots, n\}$. 
The matrix $A$ is \emph{doubly stochastic} if 
it is both row and column stochastic.  

More generally, let $\mbr = \bsmallmat r_1 \\ \vdots \\ r_m \esmallmat \in \Rm$ 
and $\mbc  = \bsmallmat c_1 \\ \vdots \\ c_n \esmallmat \in \Rn$ 
be positive vectors, that is, $r_i > 0$  for all $i \in \{1,\ldots, m\}$  
and $c_j > 0$   for all $ j\in \{1,\ldots, n\}$. 
The  $m\times n$ matrix $A = \bmat a_{i,j} \emat$ is 
\emph{$\mbr$-row stochastic} if  $A$ is  nonnegative  and 
\[
\rowsum_i(A) =  \sum_{j=1}^n a_{i,j} = r_i \qquad\text{for all $i \in \{1,\ldots, m\}$}
\]
or, equivalently, if 
\[
A \mbj_n = \mbr.
\]
The matrix $A$ is \emph{$\mbc$-column stochastic} if  $A$ is  nonnegative  and 
\[
\colsum_j(A) = \sum_{i=1}^m a_{i,j} = c_j \qquad\text{for all $j \in \{1,\ldots, n\}$}
\]
or, equivalently, if 
\[
A^t \mbj_m  = \mbc.
\]
The matrix $A$ is \emph{$(\mbr,\mbc)$-doubly stochastic} if it is both $\mbr$-row stochastic 
and $\mbc$-column stochastic.  

For positive vectors $\mbr \in \Rm$ and $\mbc \in \Rn$, 
there exists an $(\mbr,\mbc)$-doubly stochastic $m \times n$ matrix $A$ if and only if 
$\sum_{i=1}^m r_i = \sum_{j=1}^n c_j$.

An $n \times n$ matrix $\bmat x_{i,j}\emat$ is \emph{diagonal} if $x_{i,j} = 0$ for $i \neq j$.
We denote this matrix by $\diag ( x_{1,1}, \ldots, x_{n,n} ) $.  
Let $A$ be a nonnegative $m \times n$ matrix, and let $\mbc \in \Rn$  and $\mbr \in \Rm$ be positive vectors.   
If $A$ has positive column sums, then  
\[
\diag(A,\mbc) = \diag\left(\frac{c_1}{\colsum_1(A)}, \ldots, \frac{c_n}{\colsum_n(A)} \right) 
\]
is the unique $n \times n$ diagonal matrix such that 
\[
\mcc(A) =  A\diag(A,\mbc)
\]
is $\mbc$-column stochastic.  The operation $A \mapsto \mcc(A)$ is called  \emph{column scaling}.
Note that $A$ is \mbc-column stochastic if and only if $\colsum_j(A) = c_j$ for all $j \in \{1,\ldots, n\}$ 
if and only if $\diag(A,\mbc) = \diag(1,\ldots, 1) = I_n$, where $I_n$ is the $n \times n$ identity matrix, 
if and only if $\mcc(A) = A$. 
Also, if  $A$ has positive row sums, then $\mcc(A) $ has positive row sums.  

If the nonnegative $m \times n$ matrix $A$ has positive rowsums, then 
\[
\diag(A,\mbr) = \diag\left(\frac{r_1}{\rowsum_1(A)}, \ldots, \frac{r_m}{\rowsum_m(A)} \right) 
\]
is the unique diagonal matrix such that 
\[
\mcr(A) =  \diag(A,\mbr)A
\]
 is $\mbr$-row stochastic.  The  operation $A \mapsto \mcr(A)$  is called \emph{row scaling}.
The matrix $A$ is \mbr-column stochastic if and only if $\mcr(A) = A$ if and only if $\diag(A,\mbr) =I_m$.
Also, if $A$ has positive column sums, then $\mcr(A) $ has positive column sums.  

Let $\mbc \in \Rn$  and $\mbr \in \Rm$ be positive vectors.  
If the nonnegative $m \times n$ matrix $A$ has positive row and column sums, 
then we can construct an infinite sequence of matrices by alternate column and row scaling:
\beq               \label{SinkhornFinite:2}
A \mapsto \mcc(A) \mapsto \mcr(\mcc(A)) \mapsto \mcc(\mcr(\mcc(A))) \mapsto \ldots.
\eeq  
Sinkhorn and Knopp~\cite{sink-knop67} proved if $A = \bmat a_{i,j}\emat$ is a nonnegative $m\times n$ matrix 
with $\sigma$-diagonal $\prod_{i=1}^n a_{i,\sigma(i)} \neq 0$ for some permutation $\sigma \in S_n$, 
then this sequence of matrices converges coordinate-wise 
to an $(\mbr,\mbc)$-doubly stochastic matrix $S(A)$.  
Nathanson~\cite{nath2019-182} computed explicit doubly stochastic Sinkhorn limits for 
certain symmetric $3\times 3$ matrices.

For every $m \times n$ matrix $A$ with positive column and row sums and with transpose $A^t$, we have 
\[
\mcc(A)^t = \mcr\left(A^t\right) 
\qqand
\mcr(A)^t = \mcc\left(A^t\right). 
\]
These transpose relations imply that the following diagram commutes:
\[
\xymatrix{
A_0 \ar[r]^{\mcc} \ar@/_/[d]_t       &  A_1 \ar[r]^{\mcr} \ar@/_/[d]_t      & A_2  \ar[r]^{\mcc} \ar@/_/[d]_t 
& A_3 \ar[r]^{\mcr} \ar@/_/[d]_t    & \cdots \\
A_0^t \ar[r]^{\mcr}    \ar@/_/[u]_t  & A_1^t \ar[r]^{\mcc}   \ar@/_/[u]_t  & A_2^t   \ar[r]^{\mcr}   \ar@/_/[u]_t 
& A_3 \ar[r]^{\mcc}   \ar@/_/[u]_t & \cdots
}
\]
It follows that if $A$ becomes $(\mbr,\mbc)$-doubly stochastic after exactly $L$ scalings, beginning with a column scaling, 
then $A^t$ becomes $(\mbc,\mbr)$-doubly stochastic after exactly $L$ scalings, beginning with a row scaling.

There exist nonnegative matrices that become doubly stochastic after exactly two scaling operations.   
Here is an example of a positive $2\times 2$ matrix that becomes doubly stochastic 
after one row scaling and one column scaling:
\[
A = \bmat 3 &  6 \\ 5 & 10 \emat \longrightarrow  \bmat 1/3 &  2/3 \\ 1/3 & 2/3 \emat \longrightarrow 
 \bmat 1/2 &  1/2 \\ 1/2 & 1/2 \emat.  
\]
One column scaling followed by one row scaling also produces a doubly stochastic matrix:
\[
A = \bmat 3 &  6 \\ 5 & 10 \emat \longrightarrow  \bmat 3/8 &  3/8 \\ 5/8 & 5/8 \emat \longrightarrow 
 \bmat 1/2 &  1/2 \\ 1/2 & 1/2 \emat. 
\]
Nathanson~\cite{nath2020-186} proved that if a nonnegative $2\times 2$ matrix becomes doubly stochastic 
after a finite number of scalings, then it becomes doubly stochastic after at most two scalings, 
and determined all $2\times 2$  matrices with this property.  
He asked if there also exist positive $3\times 3$ matrices that become doubly stochastic 
after a finite number of row and column scalings.
Ekhad and Zeilberger~\cite{ekha-zeil19} constructed a positive $3\times 3$ matrix with the property 
that it becomes doubly stochastic after exactly one row scaling and one column scaling.   
Nathanson~\cite{nath2019-184} extended this construction to obtain $n\times n$ matrices 
for all $n \geq 3$ that become doubly stochastic after two scalings.
 He asked if for every $n$ there exists an integer $L(n)$ such that, if a nonnegative $n\times n$ 
 matrix converges after only $L$ scalings, then $L \leq L(n)$.

 In this paper we prove that $L(n) \leq 2$ for all $n \geq 2$.
 Thus, there does not exist a nonnegative $n \times n$ matrix that becomes doubly stochastic after 
exactly three scaling operations.  
The matrix obtained after one scaling operation is either row or column stochastic.  
Thus, it suffices to prove that there does not exist a matrix that is column stochastic (resp. row stochastic) 
but not row stochastic (resp. column stochastic), 
and that becomes doubly stochastic after exactly two scaling operations.  
More generally, we prove this for $(\mbr,\mbc)$-doubly stochastic matrices.

\section{Finitely many scalings}

We use the following elementary inequality.

\bl                \label{SinkhornFinite:lemma:z}
Let $c_1,\ldots, c_n$ be positive real numbers and let $z_1,\ldots, z_n$ be real numbers, 
not all 0,  such that $z_i + 1 > 0$ for $i \in \{1,\ldots, n\}$.  
If  $\sum_{i=1}^n c_i z_i = 0$, then  
\[
\sum_{i=1}^n \frac{c_i z_i}{z_i + 1} < 0.
\]
\el

\begin{proof}
Because  $\sum_{i=1}^n c_i z_i = 0$,  there exist $j,k \in \{1,\ldots, n\}$ such that $z_j > 0$ and $z_k < 0$.  
Renumbering, we can assume that there exists $r \in \{1,\ldots, n-1\}$ such that 
$z_i \geq 0$ for   $i \in \{1,\ldots, r \}$ and $z_i < 0$ for   $i \in \{r+1,\ldots, n\}$.  
For $i \in \{ 1,\ldots, r \}$, we have $z_i+1 \geq 1$ and so 
\[
\frac{c_i z_i}{z_i + 1} \leq c_iz_i.
\]
For $i \in \{  r+1,\ldots, n \}$, we have $0 < z_i+1 < 1$ and so 
\[
\frac{c_i |z_i|}{z_i + 1} >  c_i |z_i| \qqand \frac{c_i z_i}{z_i + 1} <  c_i z_i.
\]
Therefore, 
\[
\sum_{i=1}^n \frac{c_i z_i}{z_i + 1} < \sum_{i=1}^n c_i z_i = 0.  
\]
This completes the proof.  
\end{proof}

We also use the following result, sometimes called the
 \emph{second fundamental theorem of linear algebra}.

\bl        \label{SinkhornFinite:lemma:kernel-image}
The kernel of $A$  is the orthogonal complement of the image of $A^t$ in \Rn.
The kernel of $A^t$  is the orthogonal complement of the image of $A$ in \Rm.
\el

\begin{proof}
Strang~\cite[p. 198]{stra09}.  
\end{proof}

\bt
Let $A$ be a nonnegative $m \times n$ matrix with positive row and column sums.  
Let  $\mbr = \bsmallmat r_1 \\ \vdots \\ r_m \esmallmat \in \Rm$  
and $\mbc  = \bsmallmat c_1 \\ \vdots \\ c_n \esmallmat \in \Rn$ be positive vectors.  
If alternate column and row scaling transforms $A$ into an 
$(\mbr,\mbc)$-doubly stochastic matrix in exactly $L$ steps, then $L \leq 2$.
\et

Note that a  scaling operation is either a row scaling or a column scaling.  
A row scaling followed by a column scaling counts as two scaling operations.  

\begin{proof}
It is equivalent to prove that there does not exist a nonnegative $m \times n$ matrix $A_1$ 
with positive row and column sums that satisfies the following properties:
\benum
\item
$A_1$ is \mbc-column stochastic but not  \mbr-row stochastic, 
\item
$A_2 = \mcr(A_1)$  is \mbr-row stochastic but not \mbc-column stochastic, 
\item
$A_3 = \mcc(A_2)$ is $(\mbr,\mbc)$-doubly stochastic.
\eenum

Suppose there exists a matrix $A_1$  satisfying conditions (1), (2), and (3):  
\[
\xymatrix{
A_1 \ar[r]^{\mcr}  & A_2 \ar[r]^{\mcc} & A_3.
}
\]
Let $r_{1,i} = \rowsum_1(A_1)$ for $i \in \{1,\ldots, m\}$ 
and let $c_{2,j} = \colsum_j(A_2)$ for $j \in \{1,\ldots, n\}$. 

Note that $r_{1,i}/r_i> 0$ for all $i \in \{1,\ldots, m\}$. 
Because the matrix $A_1$ is not \mbr-row stochastic, 
there exists $i \in \{1,\ldots, m\}$ such that 
%$r_{1,i} \neq r_i$, or, equivalently, 
$r_{1,i} / r_i \neq 1$.
We have  
\[
A_2 = \mcr(A_1) =  \diag\left(r_1/r_{1,1},\ldots, r_m/r_{1,m} \right) A_1
\]
and so 
\[
A_1 = D_1 A_2 
\]
with 
\[
D_1 = \diag\left(r_{1,1}/r_1,\ldots, r_{1,m} /r_m\right) \neq I_m.
\]

Similarly, $c_{2,j} / c_j  > 0$  for all $j \in \{1,\ldots, n\}$.
Because the matrix $A_2$ is not \mbc-column stochastic, 
there exists $j \in \{1,\ldots, n\}$ such that 
$c_{2,j} /c_j  \neq 1$.  We have 
\[
A_3 = \mcc(A_2) = A_2\diag\left( c_1/c_{2,1},\ldots, c_n/c_{2,n} \right)
\]
and so 
\[
A_2= A_3 D_2 
\]
with
\beq        \label{SinkhornFinite:D2}
D_2 =  \diag\left(c_{2,1}/c_1,\ldots, c_{2,n} /c_n \right)  
\neq   I_n.
\eeq
It follows that 
\beq        \label{SinkhornFinite:DAD}
A_1 = D_1 A_3 D_2.  
\eeq

The matrices $A_2$ and $A_3$ are $\mbr$-row stochastic, and so 
\[
\mbr = A_3 \mbj_n =  A_2 \mbj_n = A_3D_2 \mbj_n.  
\]
Therefore, 
\[
A_3 \left( D_2 \mbj_n - \mbj_n \right) = \mbo 
\]
and   
\beq                              \label{SinkhornFinite:kernel}
\bmat z_1 \\ \vdots \\ z_n \emat = \mbz = D_2 \mbj_n - \mbj_n 
 = \bmat c_{2,1}/c_1  - 1\\ \vdots \\ c_{2,n}/c_n -1\emat 
 \in \kernel(A_3).  
\eeq
We have $\mbz \neq \mbo$ and 
\[
 \bmat z_1 + 1 \\ \vdots \\ z_n + 1 \emat  = \mbz + \mbj_n = D_2 \mbj_n 
 =  \bmat c_{2,1}/c_1 \\ \vdots \\ c_{2,n}/c_1 \emat > \mbo.
\]
Thus, 
\[
D_2 = \diag(z_1+1,\ldots, z_n+1)
\]
and
\[
D_2^{-1} = \diag\left( \frac{1}{z_1+1},\ldots, \frac{1}{z_n+1} \right).
\]
We obtain  
\[
D_2^{-1} \mbc = \diag\left( \frac{1}{z_1+1},\ldots, \frac{1}{z_n+1} \right)\bmat c_1\\ \vdots \\ c_n \emat 
= \bmat   \frac{c_1}{z_1+1} \\ \vdots  \\ \frac{c_n}{z_n+1}\emat. 
\]

The matrices $A_1$ and $A_3$ are $\mbc$-column stochastic, and so
\[
\mbj_m^t A_1 = \mbj_m^t A_3 = \mbc^t.
\]
Because $\mbz \in \kernel(A_3)$, we have  $A_3 \mbz = \mbo$ and so 
\beq            \label{SinkhornFinite:cz-sum}
0 = \mbj_m^t \mbo = \mbj_m^t A_3 \mbz = \mbc^t \mbz = \sum_{j =1}^n c_j z_j.
\eeq 
Because $A_1$ is $\mbc$-column stochastic, identity~\eqref{SinkhornFinite:DAD} implies 
\[
\mbc =  A_1^t \mbj_m = D_2 A_3^tD_1\mbj_m. 
\]
Equivalently, 
\beq               \label{SinkhornFinite:image}
A_3^t\left( D_1\mbj_m \right)  = D_2^{-1} \mbc = \bmat \frac{c_1}{z_1 + 1} \\ \vdots \\ \frac{c_n}{z_n + 1} \emat 
\in \image(A_3^t).
\eeq
By Lemma~\ref{SinkhornFinite:lemma:kernel-image}, 
\[
 \image(A_3^t) =  \left( \kernel(A_3)\right)^{\perp}.
\]
Relations~\eqref{SinkhornFinite:kernel} and~\eqref{SinkhornFinite:image}  imply that 
\beq            \label{SinkhornFinite:cz+1-sum}
0 = (D_2^{-1} \mbc, \mbz) = \sum_{i=1}^n \frac{c_i z_i}{z_i + 1}.
\eeq
Equations~\eqref{SinkhornFinite:cz-sum} and~\eqref{SinkhornFinite:cz+1-sum} 
contradict Lemma~\ref{SinkhornFinite:lemma:z}, 
and so  matrices $A_1$, $A_2$, and $A_3$ satisfying conditions~(i),~(ii), and~(iii) do not exist.   
This completes the proof.  
\end{proof}

\def\cprime{$'$} \def\cprime{$'$}
\providecommand{\bysame}{\leavevmode\hbox to3em{\hrulefill}\thinspace}
\providecommand{\MR}{\relax\ifhmode\unskip\space\fi MR }
\providecommand{\MRhref}[2]{%
  \href{http://www.ams.org/mathscinet-getitem?mr=#1}{#2}
}
\providecommand{\href}[2]{#2}


\begin{thebibliography}{1}


\bibitem{ekha-zeil19}
S.~B. Ekhad and D.~Zeilberger, Answers to some questions about explicit
  Sinkhorn limits posed by Mel Nathanson, {arXiv:1902.10783}, 2019.

\bibitem{nath2019-182} 
M.~B.~Nathanson, Matrix scaling  and explicit doubly stochastic limits, 
\emph{Linear Algebra and its Applications} 578 (2019), 111--132.  

\bibitem{nath2019-184} 
M.~B.~Nathanson, Matrix scaling limits in finitely many iterations, 
in:  \emph{Combinatorial and Additive Number Theory III},
Springer, New York, 2019.  

\bibitem{nath2020-186} 
M.~B.~Nathanson, Alternate minimization and doubly stochastic matrices, 
Integers 20 (2020), to appear.


\bibitem{sink-knop67}
R. Sinkhorn and P. Knopp, 
Concerning nonnegative matrices and doubly stochastic matrices, 
Pacific J. Math. 21 (1967), 343--348.


\bibitem{stra09}
{G. Strang}, \emph{{Introduction to Linear Algebra}}, 4th ed.,
  {Wellesley-Cambridge Press}, {Wellesley, MA}, 2009.

\end{thebibliography}
\end{document}